\documentclass[final,3p,times]{elsarticle}

\usepackage{amssymb}

\usepackage{lineno}

\journal{Journal of the Mechanics and Physics of Solids}
\usepackage{booktabs}
\usepackage{amsmath}
\usepackage{amsthm}
\usepackage{bm}
\usepackage{float}
\usepackage{xcolor}
\usepackage{subfigure}
\usepackage{graphicx}
\usepackage{makecell}
\usepackage{multirow} 
\usepackage{threeparttable}
\usepackage[strings]{underscore}
\usepackage[figuresright]{rotating}
\usepackage{siunitx}
\usepackage{adjustbox}
\biboptions{numbers,sort&compress}
\usepackage{cleveref} 
\PassOptionsToPackage{numbers}{natbib}
\usepackage{natbib}

\newcommand{\G}{\mathbf{G}}
\newcommand{\g}{\mathbf{g}}
\newcommand{\Gab}{G_{\alpha \beta}}
\newcommand{\gab}{g_{\alpha \beta}}

\newcommand{\psiel}{\psi}

\definecolor{jf_red}{rgb}{0, 0, 0}
\definecolor{rv}{rgb}{0, 0, 0}

\newtheoremstyle{boldremark}
{} 
{} 
{}  
{}          
{\bfseries} 
{.}         
{.5em}      
{}          
\theoremstyle{boldremark}
\newtheorem{remark}{Remark}

\begin{document}

\begin{frontmatter}

\title{A variationally consistent membrane wrinkling model based on spectral decomposition of the stress tensor}
\author{Daobo Zhang\fnref{label1}\corref{cor1}}
\ead{daobo.zhang@unibw.de}
\affiliation[label1]{
	organization={Institute of Engineering Mechanics and Structural Analysis - University of the Bundeswehr Munich},
	addressline={Werner-Heisenberg-Weg 39}, 
	postcode={85579}, 
	state={Neubiberg},
	country={Germany}
}

\cortext[cor1]{Corresponding author.}
\author[label1]{Josef Kiendl}

\begin{abstract}
We present a variationally consistent wrinkling model based on spectral decomposition of the stress tensor, providing a unified formulation that captures the three distinct membrane states. Compared to the previous strain-based spectral decomposition approach~\cite{zhang2024variationally}, the proposed model improves accuracy by satisfying the uniaxial tension condition from tension field theory and aligning with the mixed wrinkling criterion. It also demonstrates excellent performance under various loading conditions and offers enhanced generality by unifying strain-based, stress-based, and mixed criteria within a single framework. Beyond these improvements, the model retains the superior convergence properties of the previous approach, including the framework for the flexible inclusion or omission of residual compressive stiffness. This mitigates non-convergence or singularities in slackening states. With these adjustments, new expressions for stress and constitutive tensors are consistently derived. Finally, extensive validation through analytical, numerical, and experimental benchmark tests highlights the robustness of the model. The results confirm its accuracy in capturing the mechanical response of wrinkled thin membranes, strong convergence properties, and value for advanced membrane wrinkling analysis.

\end{abstract}

\begin{highlights}
\item Research highlight 1 (A variationally consistent wrinkling model based on stress decomposition)
\item Research highlight 2 (The model unifies strain-, stress-, and mixed-based criteria within a single framework)
\item Research highlight 3 (It ensures physical correctness while integrating the strengths of different approaches)
\item Research highlight 4 (Residual compressive stiffness enables slackened states and improves numerical stability)
\item Research highlight 5 (Simple implementation and very good iterative convergence)
\end{highlights}

\begin{keyword}
Numerical method; Stress decomposition; Membrane; Wrinkling
\end{keyword}

\end{frontmatter}


\clearpage
\section{Introduction}
\label{intro}

Wrinkling in thin membrane structures refers to out-of-plane displacements caused by localized instability under compressive stress. This phenomenon is common in elastic membranes and affects their mechanical behavior, stability, and functionality. It plays a crucial role in various engineering applications, including aerospace (e.g., solar sails, parachutes), architecture, flexible electronics, and biomedical devices. Despite extensive research, modeling and predicting wrinkling remain challenging.

To predict the mechanical behavior of wrinkled membranes, various analytical solutions have been developed. These analytical approaches can be mainly divided into two categories: those based on tension field theory (TFT), as described in \cite{wagner1931flat, reissner1938tension, Mansfield1969, Wu1978,stein1961analysis,mikulas1964behavior}, and those based on F{\"o}ppl-von K{\'a}rm{\'a}n theory, as detailed in \cite{Cerda2003, WesleyWong2006a, Puntel2011,friedl2000buckling,jacques2005mode}. The former assumes that the stress field is uniaxial tensile when wrinkling occurs in a membrane, meaning that the first principal stress is positive while the second principal stress is zero. These solutions are primarily used to describe stress fields in wrinkled membranes under small deformation problems. In contrast, the latter can resolve the detailed characteristics of wrinkles, such as their wavelength and amplitude. This approach has garnered significant attention, particularly for stretching-induced wrinkling problems. Moreover, extensions of these methods provide effective theoretical guidance for determining critical tensile stress, critical values of applied stretch, and other vital parameters. A detailed review can be found in~\citep{wang2022mechanics}. However, the aforementioned analytical solutions are generally limited to specific geometries and boundary conditions, making them difficult to extend to broader applications. Therefore, numerical methods are preferred when dealing with problems involving complex geometries and boundary conditions due to their flexibility and adaptability.

Numerical analyses of membrane wrinkling problems typically fall into two types, depending on the required detail in capturing wrinkle characteristics. The first type uses nonlinear buckling analysis with shell elements to explicitly model wrinkle morphology, accurately capturing features such as wavelength and amplitude~\citep{WesleyWong2006b,Taylor2014,Verhelst2021}. However, this approach demands extremely fine meshes, resulting in high computational costs and challenges due to scale mismatches between wrinkles and mesh resolution~\cite{Jarasjarungkiat2008}. It also often requires introducing initial geometric imperfections, adding complexity and uncertainty~\cite{Iwasa2004}. Additionally, wrinkling can be also predicted with other numerical models based on the nonlinear Föppl-von Kármán (FvK) thin plate theory~\cite{kim2012numerical,fu2019modeling,fu2021computing} and Koiter’s nonlinear plate theory ~\cite{steigmann2013well,taylor2014spatial,taylor2015comparative}. However, these methods still face high costs for finer meshes when handling short wavelengths. The Fourier reduced model~\cite{damil2010influence,damil2013new,huang2019fourier}, which discretizes the FvK equations via a Fourier series, efficiently predicts macroscopic deformations and wrinkling patterns with fewer degrees of freedom, however, its broader applicability requires further validation.

\begin{figure}[htbp]
	\centering
	\includegraphics[scale=.60]{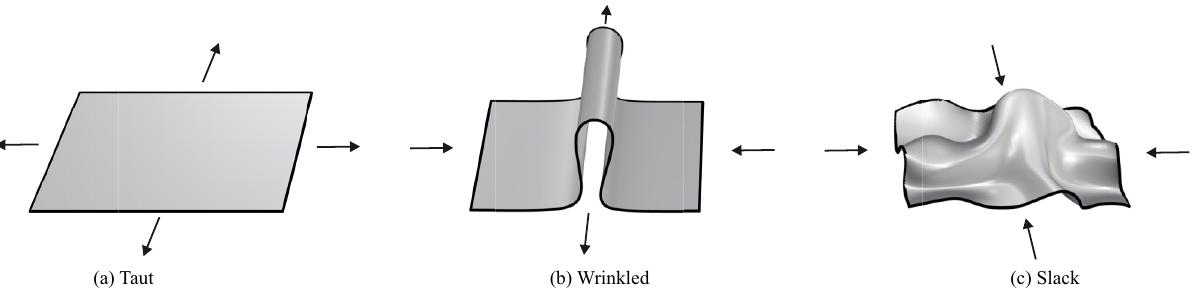}
	\caption{Possible states of a membrane: (a) taut, (b) wrinkled, (c) slack.}
	\label{FIG:states}
\end{figure}

The second type is rooted in tension field theory, which defines the stress in the wrinkle direction as zero when a membrane undergoes wrinkling. It employs coarse meshes of membrane elements with appropriate wrinkling models. Rather than explicitly resolving each wrinkle's precise geometry, this approach focuses on identifying wrinkle regions and orientations while effectively predicting the evolution of overall stress and displacement fields. Although it is less capable of capturing detailed wrinkle features, its significantly improved computational efficiency makes it more suitable for large-scale and geometrically complex problems.  According to the principle of variational consistency, existing wrinkling models can be broadly classified into two subtypes. The first subtype comprises non-variational consistent models, which explicitly rely on wrinkling criteria. These models first classify membrane states into three categories—taut, wrinkled, and slack—based on the chosen wrinkling criterion, as illustrated in the Fig.~\ref{FIG:states}. For each membrane state, distinct stiffness expressions are provided. Depending on the wrinkling criteria, these models can be further divided into those dominated by principal stress~\citep{Contri1988, nishima1986membrane, tabarrok1992nonlinear}, principal strain~\citep{Miller,Miller1982}, or a hybrid of principal stress and strain~\citep{Liu2001,Rossi2005,Feng2014}. The mixed criterion is particularly favored, as it provides a physically more accurate assessment of the membrane state compared to the other two criteria~\citep{Nakashino2005}.

Furthermore, these models can also be subdivided into kinematic modification methods and material modification methods. The kinematic modification methods focus on introducing additional parameters, such as wrinkle angle and wrinkle intensity, to adjust the deformation gradient tensor. This adjustment results in a modified strain tensor that can be decomposed into elastic strain and wrinkle strain, ensuring that the stress field satisfies the uniaxial tension condition and eliminating compressive stress in the wrinkle direction. For instance, the model presented in ~\citep{Roddeman1987a, Roddeman1987b, Raible2005, Schoop2002, Lu2001} achieves this by incorporating virtual elongation into the deformation gradient tensor, while the model introduced in ~\citep{Miyazaki2006} further accounts for shear deformation in addition to virtual elongation. Additionally, the approach developed in ~\citep{Kang1997, Kang1999} directly modifies the Green-Lagrange strain tensor, eliminating the need for complex geometric parameterization typically required by other models. In comparison, material modification methods focus more on refining the constitutive law to minimize or eliminated compressive stresses along the wrinkle direction. A representative early approach is the iterative membrane properties model (IMP), proposed by \cite{Miller,Miller1982}. This method iteratively adjusts the local elasticity matrix based on the previous load increment, allowing it to adapt to the stress state in wrinkled regions. In addition to the IMP, other approaches include penalization techniques and projection methods. Penalization techniques introduce a penalty parameter to adjust the relevant components of the material tensor~\citep{Contri1988,Liu2001,Woo2004,Rossi2005,Jarasjarungkiat2008}, while projection methods construct a projection matrix that interacts with the material tensor~\citep{Akita2007,Jarasjarungkiat2009,LeMeitour2021}. Both approaches work by softening the stiffness of the material in the compressive direction, effectively reducing it to near zero in the wrinkle direction.

The second subtype encompasses variationally consistent models, primarily represented by Pipkin’s model~\cite{Pipkin1986,pipkin1993convexity,pipkin1994relaxed}. Within a variational framework, tension field theory was reformulated as an energy minimization problem to describe membrane wrinkling. The strain energy function was assumed to be quasiconvex, and replaced by a relaxed energy functional, ensuring that no compressive stresses appear in the solution. Building on this foundation,~\citep{Mosler2008,Mosler2009} further refined the model and implemented it in a finite element code, resulting in a computationally efficient formulation. Recently,~\citep{zhang2024variationally} proposed a novel wrinkling model that achieves variational consistency by spectrally decomposing the strain tensor to split the strain energy into tensile and compressive parts. This approach employs a unified expression that satisfies the stiffness requirements for different membrane states and demonstrates promising predictive accuracy and convergence performance in various benchmark tests. However, as this model relies on a wrinkling criterion based on principal strain, its stress predictions are still influenced by Poisson’s effect. Therefore, it is essential to develop a new variational consistent wrinkling model that builds upon the mixed criterion widely acknowledged as the most effective.

In this paper, we present a novel variationally consistent wrinkling model, which utilizes spectral decomposition of the stress tensor and adheres to the mixed wrinkling criterion. 
The paper is structured as follows: Section \ref{formulation} outlines the mechanical fundamentals of membrane structures, providing the theoretical basis for subsequent sections and briefly reviews the underlying theory and definitions of existing wrinkling models, while Section \ref{new_wrinkiling_model} introduces the new wrinkling model. Section \ref{IGA} summarizes the linearization of the variational formulation for the model and the discretization. To demonstrate the applicability of our approach, Section \ref{results} presents several examples, including benchmark cases with analytical, experimental, or numerical reference solutions. Finally, Section \ref{conclusion} highlights the main characteristics of the proposed wrinkling model and discusses potential directions for future research.
\section{Membrane formulation and theoretical background of wrinkling models}
\label{formulation}
\begin{figure}[H]
	\centering
	\includegraphics[scale=.70]{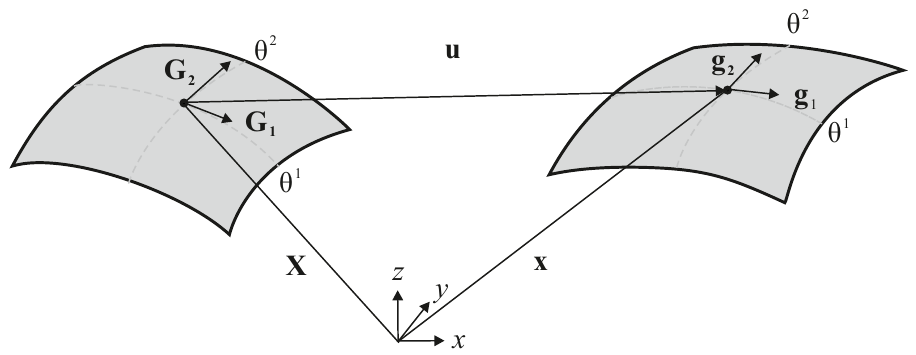}
	\caption{Schematic view of a membrane \cite{zhang2024variationally}. The vectors $\mathbf{X}$ and  $\mathbf{x}$ denote the positions on the midsurface in the reference and deformed configurations, respectively, with $\mathbf{u}$ as the displacement.} 
	\label{FIG:geometry}
\end{figure}

This section outlines the key aspects of membrane geometry, kinematics, and constitutive law. Additionally, it provides a brief introduction to the theoretical background utilized in existing wrinkling models. As illustrated in Fig.~\ref{FIG:geometry}, a membrane is typically modeled as a curved surface with thickness 
$t$, described by the position vector $\mathbf{x}\left(\theta^{\alpha}\right)$. The parameters $\theta^{\alpha}$ denote the surface coordinates, where Greek indices run over $\left\lbrace 1,2\right\rbrace $, and repeated indices follow the Einstein summation convention. In the current configuration, the covariant base vectors $\g_{\alpha}$ on the midsurface are defined by:
\begin{equation}
	\g_{\alpha} = \frac{\partial \mathbf{x}}{\partial \theta^{\alpha}}=\mathbf{x}_{,\alpha}. 
\end{equation}
In the reference configuration, these covariant base vectors are denoted $\mathbf{G}_{\alpha}$. The corresponding contravariant vectors, $\G^{\alpha}$ and $\g^{\alpha}$, satisfy
\begin{equation}
	\G^{\alpha} \cdot \G_{\beta}= \delta^{\alpha}_{\beta}, \quad  \g^{\alpha} \cdot \g_{\beta}= \delta^{\alpha}_{\beta},
\end{equation}
where $\delta^{\alpha}_{\beta}$ is the Kronecker delta, equal to $1$ if $\alpha = \beta$ and $0$ otherwise. The metric coefficients in the reference and current configurations are then expressed as:
\begin{equation}
	\Gab = \G_{\alpha}\cdot \G_{\beta}, \quad  \gab = \g_{\alpha}\cdot \g_{\beta}.
\end{equation}
The deformation gradient $\mathbf{F}$ is subsequently introduced in the form:
\begin{equation}
	\label{eq_F}
	\mathbf{F} =  \g_{\alpha} \otimes \G^{\alpha}, \quad \mathbf{F}^T =  \G_{\alpha} \otimes \g^{\alpha}.
\end{equation}
To handle nonlinear strain-displacement relationships, the Green--Lagrange strain tensor $\mathbf{E}$ is used:
\begin{equation}
\mathbf{E}= \frac{1}{2}\left(\mathbf{F}^T\mathbf{F}- \mathbf{I}\right),
\end{equation}
with $\mathbf{I}$ as the identity tensor. In particular, its in-plane components take the form:
\begin{equation}
	E_{\alpha\beta} =\frac{1}{2}\left(\gab- \Gab\right).
\end{equation}
The second Piola–Kirchhoff (PK2) stress tensor $\mathbf{S}$, which is energetically conjugate to $\mathbf{E}$, is derived from the strain energy density $\psiel$. Under plane-stress conditions for a St.~Venant–Kirchhoff material, $\psiel$ takes the form:
\begin{equation}
	\label{str_energy_org}
	\psiel\left(\mathbf{E}\right) = \frac{\lambda}{2} \left(\text{tr}\left(\mathbf{E}\right)\right)^{2}+\mu \text{tr}\left(\mathbf{E}^{2}\right) - \frac{\lambda^2}{2\left(\lambda+2\mu\right)} \left(\text{tr}\left(\mathbf{E}\right)\right)^{2},
\end{equation}
where $\lambda$ and $\mu$ are the Lamé parameters. The corresponding plane-stress material tensor $\mathbb{C}^{\alpha\beta\gamma\delta}$ provides the stress components:
\begin{equation}
	S^{\alpha\beta} = \mathbb{C}^{\alpha\beta\gamma\delta} E_{\gamma\delta}.
	\label{stress_org}
\end{equation}
Approximating the differential volume $dV$ by $t$ $dA$, where $dA$ is the membrane differential area, leads to the principle of virtual work:
\begin{equation}
	\begin{aligned}
		\delta W\left(\mathbf{u}, \delta \mathbf{u}\right)  =\delta W^{\mathrm{int}}-\delta W^{\mathrm{ext}} 
		=\int_{A} \mathbf{S}: \delta \mathbf{E} \ t\ \mathrm{d}A-\int_{A} \mathbf{f} \cdot \delta \mathbf{u} \ \mathrm{d}A,
	\end{aligned}
	\label{weak form}
\end{equation}
Here, $\delta \mathbf{E}$ is the virtual strain, $\mathbf{f}$ represents external forces, and $\delta \mathbf{u}$ is the virtual displacement. The two integrals account for the internal work $\delta W^{\text{int}}$ and external work $\delta W^{\text{ext}}$, respectively. Further discretization can be found in Section~\ref{IGA}.

Next, we briefly introduce the theoretical background of existing wrinkling models. These models are based on three main components: (1) the wrinkling criterion, which classifies the membrane state as taut, wrinkled, or slack; (2) the uniaxial tension condition, derived from tension-field theory, assuming zero stress along the wrinkling direction; and (3) the determination of the wrinkling direction, typically formulated as a nonlinear equation. However, since this study is confined to the isotropic material framework, the wrinkling direction naturally aligns with the eigenvector of the stress or strain tensor. Therefore, this aspect is not further discussed. For more details, see \citep{Miyazaki2006, Jarasjarungkiat2008, Nakashino2005}.
\begin{figure}[H]
	\centering
	\includegraphics[scale=.75]{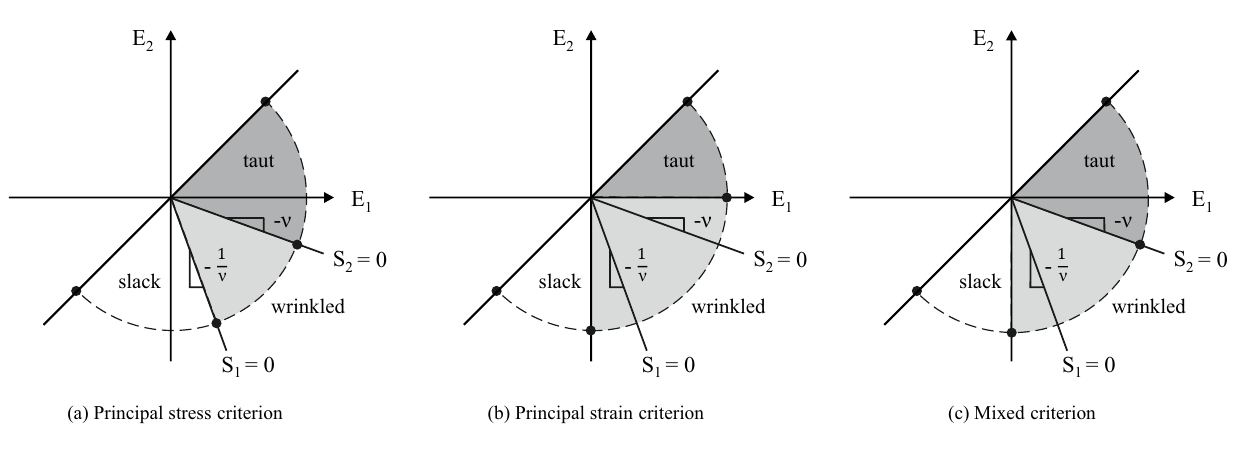}
	\caption{Comparison of membrane state classifications by wrinkling criteria.} 
	\label{FIG:criteria}
\end{figure}

As shown in Fig.~\ref{FIG:criteria}, wrinkling criteria define the membrane state under various loading conditions and are categorized based on principal stress, principal strain, or a combination of both. The principal stresses and principal strains are denoted as $S_{1/2}$ and $E_{1/2}$, respectively, with the assumption that $S_1 \geq S_2$ and $E_1 \geq E_2$. In some cases, the principal stress and strain criteria may lead to misclassifications due to Poisson’s effect, while the mixed criterion provides a more accurate assessment. For instance, when the first principal stress is non-positive while the first principal strain remains positive, the principal stress criterion classifies the membrane as slack, whereas the other two criteria classify it as wrinkled. Similarly, for the taut classification, if the second principal stress is positive while the second principal strain is non-positive, the principal strain criterion identifies the membrane as wrinkled, whereas the other two criteria classify it as taut. These discrepancies arise primarily due to the effects of Poisson’s ratio. Consequently, wrinkling models based on different wrinkling criteria may produce correspondingly varied results. For further discussion, see \citep{Feng2014, Rossi2005}.

\begin{figure}[htbp]
	\centering
	\includegraphics[scale=.70]{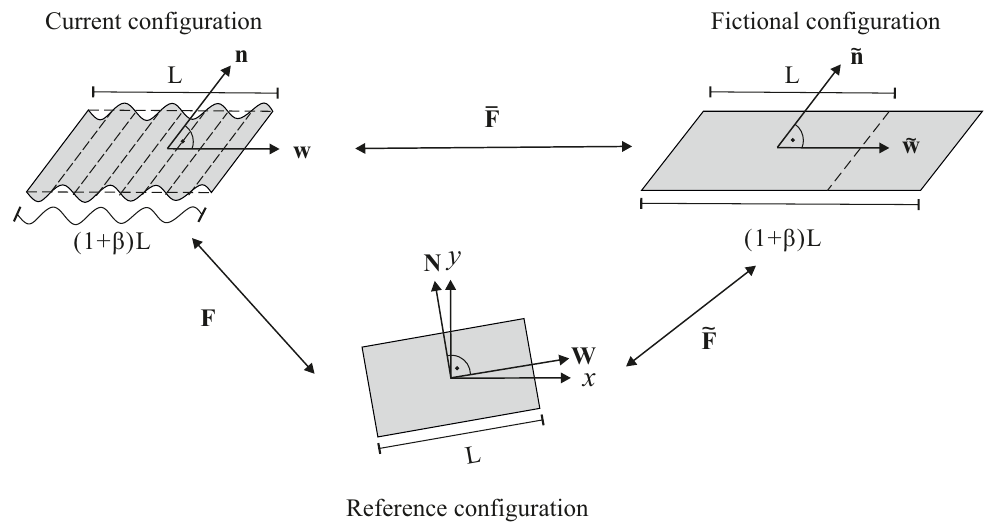}
	\caption{Transformations among the various configurations of a membrane of length $L$.} 
	\label{FIG:wrinkle_configuration}
\end{figure}

As for the previously mentioned uniaxial tension condition, it can be expressed as:
\begin{equation}
	\label{utc}
	\tilde{\mathbf{\sigma}} \tilde{\mathbf{w}}=\mathbf{0},
\end{equation}
where $\tilde{\mathbf{\sigma}}$ is the modified Cauchy stress tensor after wrinkling and $\tilde{\mathbf{w}}$ is a unit vector that indicates the direction of the wrinkle. The Eq.~\eqref{utc} is equivalent to
\begin{equation}
	\label{eq_utc}
	\tilde{\mathbf{w}} \cdot \tilde{\mathbf{\sigma}} \tilde{\mathbf{w}}=0, \quad \tilde{\mathbf{n}} \cdot \tilde{\mathbf{\sigma}} \tilde{\mathbf{w}}=0, 
\end{equation}
where $\tilde{\mathbf{n}}$ is also a unit vector perpendicular to $\tilde{\mathbf{w}}$. Then, the Eq.~\eqref{eq_utc} can be rewritten as
\begin{equation}
	\label{eq_utc_S}
	\mathbf{W} \cdot \tilde{\mathbf{S}} \mathbf{W} =0, \quad \mathbf{N} \cdot \tilde{\mathbf{S}} \mathbf{W}=0,
\end{equation}
with the following relations:
\begin{equation}
	\mathbf{W} = \mathbf{F}^T\mathbf{w} = \mathbf{\tilde{F}}^T\mathbf{\tilde{w}}, \quad
	\mathbf{\tilde{F}} = \mathbf{\bar{F}}  \mathbf{F} = (\mathbf{I}+ \beta \mathbf{w} \otimes \mathbf{w}) \mathbf{F}, \quad	\label{stress_definition}
	\tilde{\mathbf{\sigma}}=\frac{1}{\operatorname{det} \tilde{\mathbf{F}}} \tilde{\mathbf{F}} \tilde{\mathbf{S}} \tilde{\mathbf{F}}^{\mathrm{T}},
\end{equation}
where $\mathbf{\tilde{S}}$ is the modified PK2 stress tensor in the fictional configuration, $\mathbf{F}$ the standard deformation gradient tensor in the reference configuration, $\mathbf{\tilde{F}}$ the modified deformation gradient tensor in the fictional configuration. The same applies to other quantities, as shown in Fig.~\ref{FIG:wrinkle_configuration}. The scalar $\beta$ represents the elongation factor required to stretch a wrinkled surface of length $L$ into a fictitious flat plane.  

Furthermore, the Eq.~\eqref{eq_utc_S} can be reformulated as
\begin{equation}
	\label{US}
	\hat{\tilde{S}}'^{22} = \mathbf{U}_2 \cdot \hat{\tilde{\mathbf{S}}}  = 0, \quad \hat{\tilde{S}}'^{12} = \mathbf{U}_3 \cdot \hat{\tilde{\mathbf{S}}} = 0,
\end{equation}
where $\hat{\tilde{\mathbf{S}}} = [\tilde{S}^{11}, \tilde{S}^{22}, \tilde{S}^{12}]^T$ represents the modified stress tensor $\tilde{\mathbf{S}}$ in Voigt notation, $\hat{\tilde{\mathbf{S}}}' = [\hat{\tilde{S}}'^{11}, \hat{\tilde{S}}'^{22},\hat{\tilde{S}}'^{12}]^T$ is the resulting stress in the wrinkle directions, $\mathbf{U}_2$ and $\mathbf{U}_3$ are the components of a rotation matrix related to the wrinkle directions. For further details, refer to~\cite{Miyazaki2006, Jarasjarungkiat2009, Feng2014}. 
This implies that both the normal and shear stresses in the wrinkle direction disappear, corresponding to the eigenvector direction of the stress tensor for isotropic materials. This observation leads naturally to a novel wrinkling model based on spectral decomposition of the stress tensor, which inherently yields eigenvectors well-suited for defining wrinkle directions. This approach significantly simplifies model complexity, improves computational efficiency, and enhances convergence performance. The formulation and implementation of this model are discussed in detail in the following sections.

\section{A new variationally consistent wrinkling model based on spectral decomposition of the stress tensor}
\label{new_wrinkiling_model}
In this section, we first introduce the theoretical framework of spectral decomposition \cite{Miehe2001}, which then serves as the basis for deriving our proposed wrinkling model. According to spectral decomposition, the strain tensor $\mathbf{E}$ and stress tensor $\mathbf{S}$ can be represented as:
\begin{equation}
	\label{split}
	\mathbf{E} =  E_{\alpha}\mathbf{M}^{\mathbf{E}}_{\alpha}, \quad \mathbf{S} =  S_{\alpha}\mathbf{M}^{\mathbf{S}}_{\alpha},
\end{equation}
where  $E_{\alpha}$ and ${S}_{\alpha}$ are their principal strains and stresses, respectively, and assuming $E_1 \geq E_2$ and $S_1 \geq S_2$. $\mathbf{M}^{\mathbf{E}}_{\alpha}$ and $\mathbf{M}^{\mathbf{S}}_{\alpha}$ are the second order tensors spanned with the corresponding eigenvectors $\mathbf{N}^{\mathbf{E}}_{\alpha}$ of strain tensor $\mathbf{E}$, and the eigenvectors $\mathbf{N}^{\mathbf{S}}_{\alpha}$ of stress tensor $\mathbf{S}$, respectively. They are defined respectively as
\begin{equation}
	\label{M_eig}
	\mathbf{M}^{\mathbf{E}}_{\alpha} = \frac{\partial E_{\alpha}}{\partial \mathbf{E}} = \mathbf{N}^{\mathbf{E}}_{\alpha} \otimes \ \mathbf{N}^{\mathbf{E}}_{\alpha} , \quad  \mathbf{M}^{\mathbf{S}}_{\alpha} =\frac{\partial S_{\alpha}}{\partial \mathbf{S}} = \mathbf{N}^{\mathbf{S}}_{\alpha} \otimes \ \mathbf{N}^{\mathbf{S}}_{\alpha}.
\end{equation}
Note that the index summation convention does not apply in Eq.~\eqref{M_eig}. When the material is isotropic, then the following relation holds:
\begin{equation}
	\mathbf{M}_{\alpha}=	\mathbf{M}^{\mathbf{E}}_{\alpha} =  \mathbf{M}^{\mathbf{S}}_{\alpha}.
\end{equation}
Since the wrinkle directions align with the eigenvector directions for isotropic materials, $\mathbf{M}_1$ and $\mathbf{M}_2$ represent the matrix form of $\mathbf{U}_1$ and $\mathbf{U}_2$, respectively.
Then, we can rewrite the strain energy density of plane stress for isotropic material as:
\begin{equation}
	\label{psi_org}
	\psi = \frac{1}{2}\mathbf{S}:\mathbf{E}.
\end{equation}
Building on the idea presented in \citep{Zhang2020a}, the strain energy density for linear isotropic material can be decomposed as:
\begin{equation}
	\label{psi_stress_split}
	\psi = \psi^{+} + \psi^{-}, 
\end{equation}
where the positive and negative strain energy densities $\psi^{\pm}$ are defined respectively as:
\begin{equation}
	\psi^{\pm} = \frac{1}{2}(\mathbf{S}^{\pm}:\mathbf{E}).
\end{equation}
Then, the positive and negative stresses $\mathbf{S}^{\pm}$ are defined respectively as:
\begin{equation}
	\mathbf{S}^{\pm} =   H^{\pm}(S_{\alpha})S_{\alpha} \mathbf{M}_{\alpha}, 
\end{equation}
where the principal stresses ${S}_{\alpha}$ can also be expressed in terms of principal strains using the St. Venant-Kirchhoff material law as:
\begin{equation}
	\label{stress-strain}
	{S}_{\alpha} = \frac{E}{(1-\nu^2)}({E}_{\alpha}+\nu{E}_{\beta}),\quad \alpha \ne \beta, 	
\end{equation}
with Yong's modulus $E$ and Poisson's ratio $\nu$. Thus, the positive and negative stresses can be rewritten in terms of principal strains as:
\begin{equation}
	\label{stress_split}
		\mathbf{S}^{\pm} =  \frac{E}{(1-\nu^2)} H^{\pm}(E_{\alpha}+\nu E_{\beta})(E_{\alpha}+\nu E_{\beta}) \mathbf{M}_{\alpha}, \quad \alpha \ne \beta.  
\end{equation}
Consequently, the positive and negative material tensors $\mathbb{C}^{\pm}$ can be obtained by deriving the corresponding stresses with respect to strain and are given by:
\begin{equation}
	\label{CC_stress}
	\mathbf{\mathbb{C}}^{\pm} = \frac{E}{(1-\nu^2)}\bigg( H^{\pm}(E_{\alpha}+\nu E_{\beta})\Big( (\mathbf{M}_{\alpha}+\nu\mathbf{M}_{\beta}) \otimes \mathbf{M}_{\alpha} + (E_{\alpha}+\nu E_{\beta})\frac{\partial \mathbf{M}_{\alpha}}{\partial \mathbf{E}}\Big), \quad \alpha \ne \beta. 
\end{equation}
In line with the definitions provided by \citep{Miehe2001}, the following fourth-order tensors can be expressed as:
\begin{equation}
	\mathbb{Q}_{\alpha \beta} = \mathbf{M}_{\alpha} \otimes \ \mathbf{M}_{\beta}, \quad
	\mathbb{M}_{\alpha} = \frac{\partial \mathbf{M}_{\alpha}}{\partial \mathbf{E}} = \frac{1}{2(E_{\alpha}- E_{\beta})} (\mathbb{G}_{\alpha \beta}+ \mathbb{G}_{\beta \alpha}),  \quad \alpha \ne \beta,  
	\label{QMG}
\end{equation}
with $\mathbb{G}^{ijkl}_{\alpha \beta} = M^{ik}_{\alpha} M^{jl}_{\beta} + M^{il}_{\alpha} M^{jk}_{\beta} $ and please note that the index summation convention does not apply here. Consequently, the positive and negative material tensors $\mathbb{C}^{\pm}$ can be rewritten as:
\begin{equation}
	\mathbf{\mathbb{C}}^{\pm} = \frac{E}{(1-\nu^2)} H^{\pm}(E_{\alpha}+\nu E_{\beta}) \Big( \mathbb{Q}_{\alpha \alpha} + \nu \mathbb{Q}_{\alpha \beta}  + (E_{\alpha}+\nu E_{\beta})\mathbb{M}_{\alpha} \Big), \quad \alpha \ne \beta. 
	\label{CC_stress_2}
\end{equation}
Since the negative or compressive stress has to be eliminated, then the new strain energy density $\tilde{\psi}$, its consistently derived stress $\mathbf{\tilde{S}}$ and material tensor $\tilde{\mathbb{C}}$ can be modified as:
\begin{equation}
	\tilde{\psi} = \psi^+, \quad 	\mathbf{\tilde{S}} = \mathbf{S}^+, \quad  \mathbf{\tilde{\mathbb{C}}} = \mathbf{\mathbb{C}}^+.
\end{equation}
Alternatively, in line with perspectives held by other scholars \cite{Jarasjarungkiat2008,Jarasjarungkiat2009,Miyazaki2006, Iwasa2004}, even when wrinkling occurs in thin membranes, compressive stress does not completely vanish. Instead, a certain residual compressive stress is assumed to persist. This is reflected in previous formulations, where the stress components after rotation or wrinkling, $\hat{\tilde{S}}'^{22}$ and $\hat{\tilde{S}}'^{12}$ in Eq.~\eqref{US}, are no longer strictly zero but are expressed as:
\begin{equation} \label{US_2_modified} \hat{\tilde{S}}'^{22} =\mathbf{U}_2 \cdot \hat{\tilde{\mathbf{S}}} =\eta \mathbf{U}_2 \cdot \hat{\tilde{\mathbf{S}}}, \quad \hat{\tilde{S}}'^{12} = \mathbf{U}_3 \cdot \hat{\tilde{\mathbf{S}}}= \eta \mathbf{U}_3 \cdot \hat{\tilde{\mathbf{S}}},
\end{equation} 
where $\eta$ is a small scaling factor. Building upon these approaches, we adapt our formulation to account for the persistence of compressive stress in a more general framework. Specifically, we introduce a degradation factor $\eta \ll 1$, leading to the modified expressions:

\begin{equation} \tilde{\psi} = \psi^+ + \eta \psi^-, \quad \mathbf{\tilde{S}} = \mathbf{S}^+ + \eta\mathbf{S}^- , \quad \mathbf{\tilde{\mathbb{C}}} = \mathbf{\mathbb{C}}^+ + \eta\mathbf{\mathbb{C}}^-. \end{equation}
More details on this formulation can be found in \cite{zhang2024variationally}.

Clearly, the new stress tensor resulting from the decomposition based on principal stresses fully satisfies the previously derived uniaxial tension condition, as shown in Eq.~\eqref{eq_utc_S}, where the stress in the wrinkling direction becomes zero. This follows from the fact that Eq.\eqref{eq_utc_S} can be rewritten as:
\begin{equation}
	\mathbf{W} \cdot \tilde{\mathbf{S}} \mathbf{W} =\tilde{\mathbf{S}}: \mathbf{W} \otimes \mathbf{W}=0.
\end{equation}
Since the wrinkling direction $\mathbf{W}$ in isotropic materials coincides with the principal stress direction $\mathbf{N}_2$, it follows that:
\begin{equation}
	\tilde{\mathbf{S}}: \mathbf{N}_2  \otimes \mathbf{N}_2= \tilde{\mathbf{S}}: \mathbf{M}_2  =0.
\end{equation}
Therefore, the newly derived stress tensor naturally satisfies this condition in the wrinkled state, as the basis tensors obey the orthonormal property:
\begin{equation}
	 \mathbf{M}_{\alpha} : \mathbf{M}_{\beta} =  \delta_{\alpha\beta}.
\end{equation}
 In addition, due to the nature of eigenvalue decomposition, the shear stress naturally vanishes as well.
\begin{remark}
	\label{remark1}
Following the definitions in Eq.~\eqref{QMG}, the tensors $\mathbb{Q}_{\alpha \alpha}$, $\mathbb{M}_{\alpha}$, and $\mathbb{G}_{\alpha \beta}$ are inherently major symmetric fourth-order tensors, while $\mathbb{Q}_{\alpha \beta}$ is minor symmetric. However, the sum $\mathbb{Q}_{\alpha \beta} + \mathbb{Q}_{\beta \alpha}$ remains a major symmetric tensor. Based on these properties, the modified material tensor $\tilde{\mathbb{C}}$, derived from the above stress decomposition, is major symmetric in the taut and slack states, while in the wrinkled state, it becomes minor symmetric. Notably, when Poisson's ratio $\nu = 0$, the material tensor reverts to a major symmetric form. Furthermore, this decomposition approach fully depends on a wrinkling criterion based on principal stress. As previously analyzed, compared to a mixed criterion, this model may misinterpret slack states. Specifically, when both principal stresses are less than or equal to zero and the first principal strain is positive, this model may fail to provide adequate stiffness, resulting in an overly compliant response. In contrast, as shown in~\ref{appNew}, the strain decomposition approach ensures that the material tensor remains major symmetric across all membrane states. This is because the terms $\mathbb{Q}_{12}$ and $\mathbb{Q}_{21}$ always appear or disappear simultaneously depending on the sign of $\operatorname{tr}(\mathbf{E})$, and their sum is consistently major symmetric. From this analytical perspective, we can also observe from the Eq.~\eqref{S_strain} that in the wrinkled state, where $E_1 >0$ and $E_2 \leq 0$, the uniaxial tension condition is only satisfied when $\operatorname{tr}(\mathbf{E}) \leq 0$, which results in zero stress in the wrinkling direction. Otherwise,
\begin{equation}
	\label{trE}
	\mathbf{S}^{+} : \mathbf{M}_2 = \frac{E\nu}{1-\nu^2} \operatorname{tr}(\mathbf{E}) >0, \quad \text{for}  \ \operatorname{tr}(\mathbf{E}) > 0.
\end{equation}
This suggests that, in comparison with the models based on mixed criterion (which will be introduced later), this model may overestimate the stress magnitude, leading to an overly stiff response. Furthermore, when using a wrinkling criterion based on principal strain, distinguishing between taut and wrinkled states is affected by the Poisson effect. This influence may cause the estimated stress in the taut state to be lower than the value it should ideally reach in a taut state.
\end{remark}
 
Next, we present the simplest approach that not only preserves the advantage of the stress split in satisfying the uniaxial tension condition but also converts the stress criterion into the mixed criterion by appropriately modifying the Poisson’s ratio in the aforementioned wrinkling models. The modified Poisson's ratio is given by:\begin{equation} 
	\label{slack_nu0} 
	\nu^* = H^{+}(E_2+\nu E_1)\nu. 
\end{equation}
Consequently, the Poisson effect is eliminated, meaning that the Poisson's ratio becomes zero when the membrane is in wrinkled or slack states. 
Then, substituting Eq.~\eqref{slack_nu0} into Eq.~\eqref{stress_split}, the stress formulation can be reformulated as:
\begin{equation}
	\label{mixed_split}
	\mathbf{S}^{\pm} = \frac{E}{(1-\nu^{*2})} H^{\pm}(E_{\alpha}+\nu^* E_{\beta})(E_{\alpha}+\nu^*E_{\beta}) \mathbf{M}_{\alpha}, \quad \alpha \ne \beta. 
\end{equation}
Thus, the material tensor in Eq.~\eqref{CC_stress_2} becomes
\begin{equation}
	\label{CC_mixed}
		\mathbf{\mathbb{C}}^{\pm}  = \frac{E}{(1-\nu^{*2})}  H^{\pm}(E_{\alpha}+\nu^*E_{\beta}) \Big( \mathbb{Q}_{\alpha \alpha} + \nu^* \mathbb{Q}_{\alpha\beta}  + (E_{\alpha}+\nu^* E_{\beta})\mathbb{M}_{\alpha} \Big),\quad \alpha \ne \beta. 
\end{equation}
Accordingly, the strain energy density after modification takes the form
\begin{equation}
	\label{psi_pr_strain_mix}
	\psi^{\pm} = \frac{E}{2(1-\nu^{*2})}  H^{\pm}(E_{\alpha}+\nu^* E_{\beta})(E_{\alpha}+\nu^*E_{\beta})E_{\alpha},\quad \alpha \ne \beta. 
\end{equation}

From Eq.~\eqref{CC_mixed}, we observe that this adjustment, which inherits the advantages of different decomposition approaches, brings two major improvements. First, it completely resolves the asymmetry issue introduced by the stress split. When the membrane is classified as wrinkled, the Poisson’s ratio is set to zero, eliminating $\mathbb{Q}_{\alpha\beta}$ and restoring the material tensor to a major symmetric form. Second, it fully satisfies the uniaxial tension condition and aligns with the mixed wrinkling criterion, ensuring physical consistency.

\section{Discretization}
\label{IGA}
In this section, we introduce the discretization approach, which largely follows the methodology presented in \cite{zhang2024variationally}, ensuring consistency with existing models. We then provide the essential details of the linearization and discretization procedures required to solve Eq.~(\ref{weak form}). We begin by noting that this framework supports various types of basis functions, including standard Lagrange polynomials for finite element methods or Non-Uniform Rational B-Splines (NURBS) commonly used in isogeometric analysis (IGA). Due to the many benefits of IGA, it has gained considerable traction in structural analysis. Comprehensive discussions of these advantages can be found in~\citep{cottrell2009isogeometric,hughes2005isogeometric}, so they are not reiterated here. The discretized displacement field takes the form
\begin{equation}
	\mathbf{u}=\sum_a^{n_{s h}} N^a \mathbf{u}^a,
\end{equation}
where $N^a$ are the shape functions, $n_{sh}$ is the total number of these functions, and $\mathbf{u}^a$ represents the nodal displacement vector. Its components $u^a_i \ (i=1,2,3)$ correspond to the global $x$-, $y$-, and $z$-directions. To reference each nodal displacement within a global degree of freedom index $r$, we set $r=3(a-1)+i$, yielding $u_r = u^a_i$. Differentiating with respect to $u_r$ gives
\begin{equation}
	\frac{\partial \mathbf{u}}{\partial u_r}=N^a \mathbf{e}_i,
\end{equation}
where $\mathbf{e}_i$ denotes the global Cartesian base vectors. For further details, see~\cite{Kiendl2015}. Next, by taking the derivative of the internal and external virtual work terms in Eq.~(\ref{weak form}) with respect to $u_r$, we obtain the residual force vector $\mathbf{R}$:
\begin{equation}
	R_r=F_r^{\text {int }}-F_r^{\text {ext }}=\int_{A} \mathbf{S}: \frac{\partial \mathbf{ E}}{\partial u_r} \ t \ \mathrm{~d}A-\int_{A} \mathbf{f} \cdot \frac{\partial \mathbf{u}}{\partial u_r} \mathrm{~d}A,
	\label{residual force vector}
\end{equation}
where $\mathbf{F}^{\text{int}}$ and $\mathbf{F}^{\text{ext}}$ are the nodal internal and external forces, respectively. Linearizing Eq.~(\ref{residual force vector}) leads to the tangent stiffness matrix $\mathbf{K}$, which consists of internal and external contributions, $\mathbf{K}^{\text {int }}$ and $\mathbf{K}^{\text {ext }}$:
\begin{equation}
	\begin{aligned}
		K_{rs} = K_{r s}^{\text {int }} - K_{r s}^{\text {ext }} 
		=  \int_{A} \Big( \frac{\partial \mathbf{S}}{\partial u_s}: \frac{\partial \mathbf{ E}}{\partial u_r}+\mathbf{S}: \frac{\partial^2 \mathbf{ E}}{\partial u_r \partial u_s} \Big) \ t \ \mathrm{~d}A - \int_{A} \frac{\partial \mathbf{f}}{\partial u_s} \cdot \frac{\partial \mathbf{u}}{\partial u_r} \mathrm{~d}A.
	\end{aligned}
\end{equation}
For loads that do not depend on displacements, the external stiffness matrix disappears. The linearized system is then solved iteratively via the Newton-Raphson method:
\begin{equation}
	\frac{\partial W}{\partial u_r} + \frac{\partial^2 W}{\partial u_r \partial u_s} \Delta u_s = 0,
\end{equation}
where $\Delta u_s$ are the incremental displacements. Solving the system
\begin{equation}
	\mathbf{K} \Delta \mathbf{u}=-\mathbf{R},
\end{equation}
updates the displacement field in each iteration until the residual falls below a prescribed tolerance, thereby ensuring a converged solution for the incremental displacement $\Delta \mathbf{u}$.
\pagebreak
\section{Numerical examples}
\label{results}
This section aims to evaluate the accuracy of the proposed wrinkle model predictions through five distinct numerical examples. We have selected benchmarks involving different loading conditions: pure bending, simple shear, biaxial stretching, and complex deformation of wrinkled membrane structures in three-dimensional space. The results of these numerical examples are compared against analytical solutions or well-established reference data from the literature to validate the accuracy and robustness of the proposed method. Moreover, we compare three models based on different decomposition approaches and wrinkling criteria, namely the principal stress-based model (Eq.~\eqref{CC_stress_2}), the principal strain-based model (Eq.~\eqref{C_strain}), and the mixed criterion-based model (Eq.~\eqref{CC_mixed}). All numerical simulations are based on nonlinear static isogeometric analysis using NURBS basis functions, which are particularly well-suited for accurately capturing the physical behavior of structures. In addition, we briefly describe the notations used in the following examples here. The Cauchy stress tensor, is computed using Eq.~\eqref{stress_definition}, which fully accounts for the influence of the deformation gradient $\mathbf{F}$. Additionally, the first and second principal stresses of the Cauchy stress tensor are represented by $\sigma_1$ and $\sigma_2$, respectively. For evaluating the convergence of the numerical simulations, we adopted a residual norm based on the magnitude of external forces, defined as $\Vert\mathbf{R} \Vert / \Vert\mathbf{F}^{	\text{ext}} \Vert$. This convergence criterion effectively assesses the equilibrium state during the iterative process, ensuring the desired level of computational accuracy. Unless otherwise specified, a degradation factor of $\eta = 0$ is assumed throughout all numerical analyses.
\subsection{A stretched rectangular membrane wrinkled under in-plane pure bending}
\label{exp1}
In the first numerical example, we seek to validate the proposed model by analyzing the wrinkling behavior of a pre-stretched rectangular membrane subjected to planar bending. This benchmark problem was initially formulated by \citep{stein1961analysis}, and it has since become a widely adopted standard for evaluating the performance of wrinkling models, as seen in various studies such as \cite{Lu2001, Liu2013, Feng2016, Jarasjarungkiat2009, Ding2003, Miller, GwanJeong1992}. The membrane considered in this example has dimensions defined by its length $L=2$, height $H = 1$, and thickness $t=0.01$, and possesses material properties characterized by Young's modulus $E=100$ and Poisson's ratio $
\nu=0.3$. 
\begin{figure}[htbp]
	\centering
	\includegraphics[scale=.82]{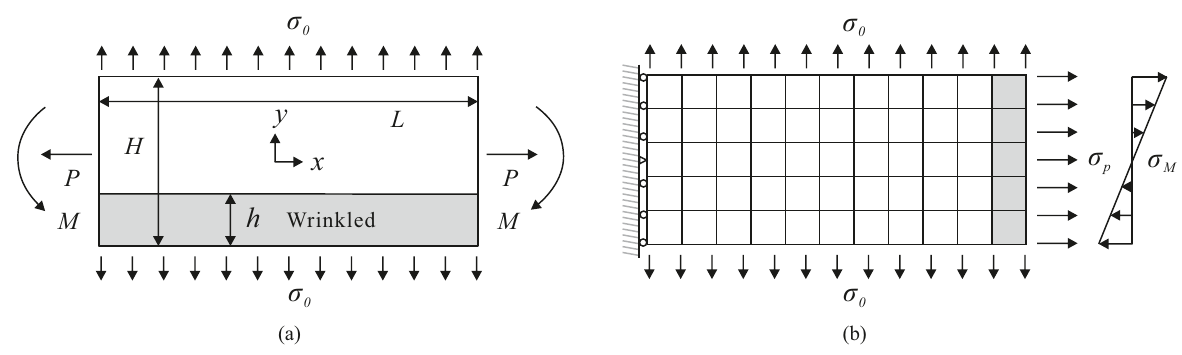}
	\caption{(a) In-plane pure bending of a pre-tensioned rectangular membrane and (b) numerical simulation settings for the right half of the rectangular membrane under in-plane bending.}
	\label{FIG:bending_settup}
\end{figure}
As illustrated in the Fig.~\ref{FIG:bending_settup}(a), the membrane experiences a uniform pre-stretching stress $\sigma_0 = 5 \times 10^{-4}$ in the $y$-direction. In addition, axial loads $P = \sigma_0 t H$ and bending moments $M$ are applied at the lateral boundaries. With the progressive increase in bending moment $M$, a wrinkling band of height $h$ starts forming along the lower edge of the membrane, as highlighted by the gray area. The height of this wrinkling band $h$ can be derived from the analytical solution provided by \cite{stein1961analysis} as:
\begin{equation}
	\frac{h}{H}= \begin{cases}0  & M / P H<1 / 6 \\ 3 M / P H-1 / 2  & 1 / 6 \leq M / P H<1 / 2\end{cases}.
\end{equation}
Due to the presence of the wrinkled region, the normal stress $\sigma_x$ is entirely eliminated in this area, while in the unwrinkled regions above, the membrane remains in a taut state. Therefore, the stress distribution along different heights can be expressed as follows:
\begin{equation}
	\frac{\sigma_{x}}{\sigma_{0}}= \begin{cases}\frac{2\left(y / H-h / H\right)}{\left(1-h / H\right)^{2}} & h / H<y / H \leq 1 \\ 0  & 0 \leq y / H \leq h / H\end{cases},
\end{equation}
where $y/H$ represents the dimensionless position within the membrane, relative to the height $H$, to describe the variation in the stress magnitude. 

Besides,  with the overall curvature $\kappa$, \citep{stein1961analysis} also provided the moment–curvature relation of a beam-like membrane, which is defined by:
\begin{equation}
	\frac{2M}{PH} = 
	\begin{cases}
		\frac{EH^2 t \kappa}{6P} \ &  \frac{EH^2 t \kappa}{2P} \leq 1 \\
		1 - \frac{2}{3} \sqrt{\frac{2P}{EH^2 t \kappa}}\  & \frac{EH^2 t \kappa}{2P} > 1
	\end{cases}.
\end{equation}

For the numerical simulation, we leveraged the symmetry of the membrane to simplify the model by analyzing only the right half, as depicted in Fig.~\ref{FIG:bending_settup}(b). Following the discretization approach suggested by \cite{Nakashino2005}, an $11 \times 5$ bi-quadratic isogeometric membrane element mesh was used. To ensure proper symmetry, the displacement in the $x$-direction along the entire left side was restricted to prevent horizontal movement, while only the midpoint of the left boundary was fixed in the $y$-direction. Equivalent stress expressions were applied to represent the axial load and bending moment, given by $\sigma_p = P/tH$ and $\sigma_{M} = 6M/tH^2 (2y/H - 1)$, respectively. Additionally, we incrementally increased the bending moment across a range of moment-load ratios $2M/PH$, varying from 0.4 to 0.825, to observe the wrinkling behavior and compare the numerical predictions to the analytical solution for model validation. It is worth noting that, in order to maintain uniform bending behavior and resist compressive stresses, the elements located at the far right end of the membrane were assumed to remain taut. Consequently, these elements were treated as standard membrane elements without incorporating wrinkling effects.

\begin{figure}[htbp]
	\centering
	\includegraphics[scale=0.9]{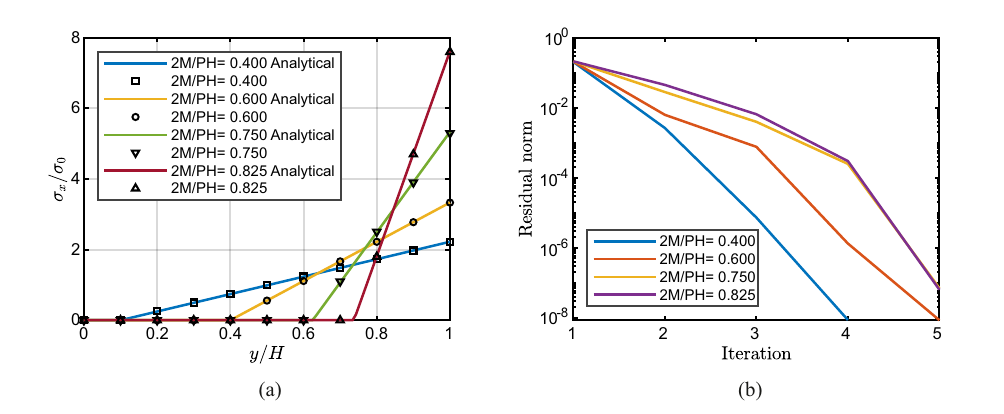}
	\caption{(a) Comparison of the normal stresses $\sigma_x$ predicted by the mixed criterion-based model with the analytical solutions under various $2M/PH$ and (b) comparison of the iteration number during the second load step.}
	\label{FIG:bending_stress}
\end{figure}

Firstly, we compared the normal stress $\sigma_{x}$ predicted by the proposed models with those derived from the analytical solution. By examining the distribution of the normal stress component $\sigma_x$, as shown in Fig.~\ref{FIG:bending_stress}(a), we used the ratio $y/H$ as the horizontal axis to present the stress response consistently across various $2M/PH$-ratios. As depicted in Fig.~\ref{FIG:bending_stress}(a), the mixed criterion-based model provided results that were consistent with the analytical solution for all tested ratios. To assess the convergence behavior of this model, the number of iterations required during the second loading step is plotted in Fig.~\ref{FIG:bending_stress}(b). This loading step exhibited the highest number of iterations among all steps. As shown in Fig.~\ref{FIG:bending_stress}(b), the mixed criterion-based model demonstrated favorable convergence properties, with a slight increase in the number of iterations needed to achieve the specified tolerance as the applied bending moment increased. Furthermore, all three models exhibit similar convergence behavior, as they are formulated within the variationally consistent framework based on spectral decomposition. Therefore, a detailed comparison of each is not necessary.

\begin{figure}[htbp]
	\centering
	\includegraphics[scale=0.9]{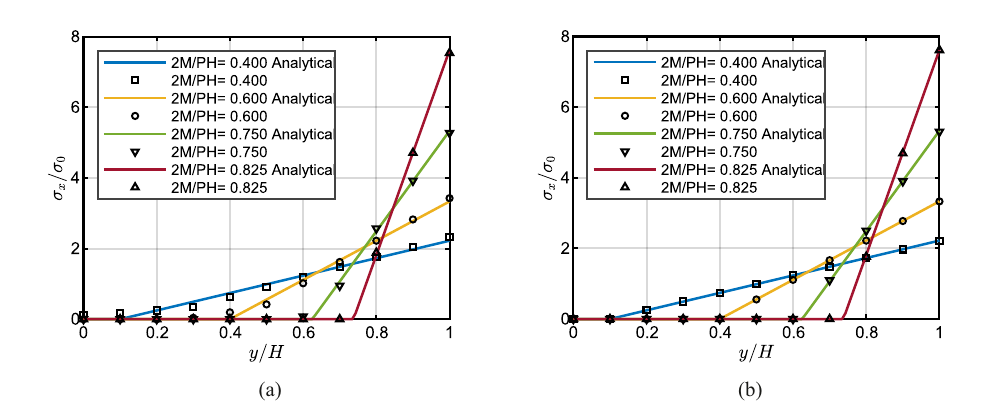}
	\caption{Comparison of the normal stresses $\sigma_x$ predicted by the models based on the strain criterion (a) and the stress criterion (b) with the analytical solutions under various $2M/PH$.}
	\label{FIG:bending_stress_strain_stress}
\end{figure}

In addition, we present the stress results predicted by the other two models. Fig.\ref{FIG:bending_stress_strain_stress}(a) shows the results of model based on the principal strain criterion, while Fig.\ref{FIG:bending_stress_strain_stress}(b) corresponds to the model based on the principal stress criterion. From these results, we observe that the predictions of the principal stress-based criterion align well with the analytical solution. In contrast, the principal strain-based criterion initially deviates from the analytical solution in the wrinkled region. Specifically, at early stages such as $2M/PH=0.4$ and $0.6$, the stress in the wrinkled region should theoretically be zero, but the model predicts a higher stress than the analytical solution. However, as the applied bending moment increases, the model gradually converges to the analytical solution. This behavior arises because the trace of the strain tensor plays a crucial role in this model, as discussed in Remark~\ref{remark1}. When the trace of the strain tensor is larger than zero, the stress in the wrinkled direction remains positive. However, when the trace of the strain tensor is zero or negative, the model satisfies the uniaxial tension condition, effectively eliminating the discrepancy with the analytical solution.

\begin{figure}[htbp]
	\centering
	\includegraphics[scale=0.9]{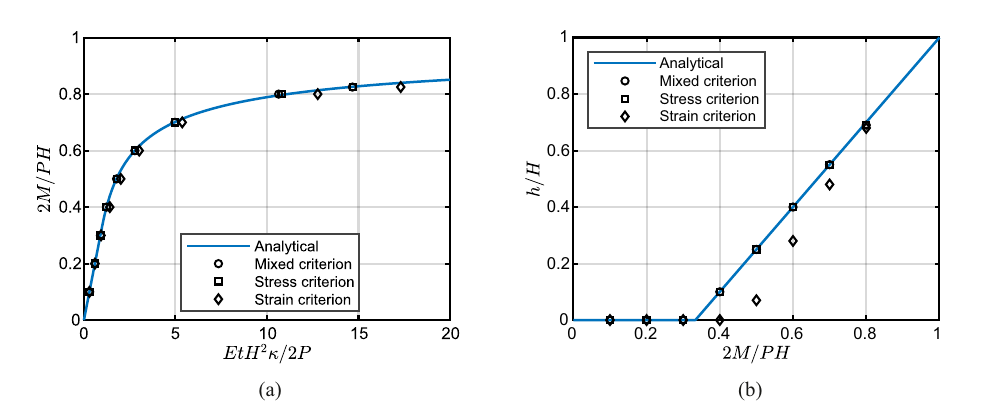}
	\caption{(a) The relation between the moment and the overall curvature and (b) the relation between the moment and the height of wrinkled band.}
	\label{FIG:bending_kappa}
\end{figure}

Furthermore, Fig.~\ref{FIG:bending_kappa} illustrates the results for the relation between the bending moment and overall curvature, see Fig.~\ref{FIG:bending_kappa}(a), as well as the relation between the bending moment and the height of the wrinkled band, see Fig.~\ref{FIG:bending_kappa}(b). These plots demonstrate that the mixed criterion-based model and the stress criterion-based model exhibit excellent agreement with the analytical solution across different bending moment-to-load ratios, both in terms of overall curvature and the predicted height of the wrinkled region. In contrast, the wrinkling model based on the strain criterion shows some deviation in curvature prediction. Regarding the height of the wrinkled band, it initially exhibits a noticeable error; however, this discrepancy gradually decreases as the bending moment increases. This is due to the same reason that causes deviations in stress prediction, as discussed earlier, and will not be reiterated here. Notably, although some deviations are present in the prediction of wrinkle band height, their impact on the stress results appears to be minimal, see Fig.~\ref{FIG:bending_stress_strain_stress}(a).

\subsection{A rectangular membrane wrinkled by simple shear}

\begin{figure}[htbp]
	\centering
	\includegraphics[scale=0.60]{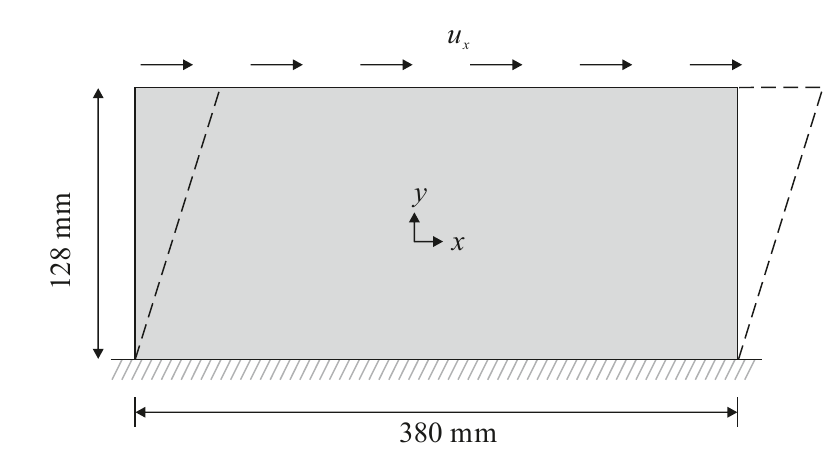}
	\caption{Geometry and boundary conditions of a rectangular membrane under shear.}
	\label{FIG:shear_test}
\end{figure}

In the second example, we evaluate the performance of the proposed models in a simple shear test applied to a rectangular membrane. This test allows us to assess the ability of the models to capture wrinkling behavior under shear deformation. The membrane has a width of $\SI{380}{\milli \meter}$, a height of $\SI{128}{\milli \meter}$ and a thickness of $\SI{0.025}{\milli \meter}$, as shown in Fig.~\ref{FIG:shear_test}. The material properties are given as Young's modulus $E = \SI{3500}{\newton \per \milli\meter\squared}$ and Poisson's ratio $\nu = 0.31$. The bottom edge is fully fixed, preventing any displacement in all directions. Additionally, the displacement in the $z$-direction is also constrained. The loading procedure follows the approach described in \cite{WesleyWong2006b,Jarasjarungkiat2009,Feng2014}: first, a vertical tensile displacement of $\SI{0.05}{\milli \meter}$ is applied to introduce the initial stiffness. Afterwards, a horizontal displacement $u_x$ of either $\SI{1.6}{\milli \meter}$ or $\SI{3}{\milli \meter}$ is incrementally applied to the top edge to induce shear deformation. Meanwhile, the vertical displacement on the top edge remains fixed, ensuring that the pre-applied stretch in the $y$-direction is preserved throughout the shear loading process. To enforce a pure shear condition, the left and right vertical edges are free. In this example, the simulation was conducted using a $40 \times 10$ mesh of bi-quadratic elements.\begin{figure}[H]
	\centering
	\includegraphics[scale=0.85]{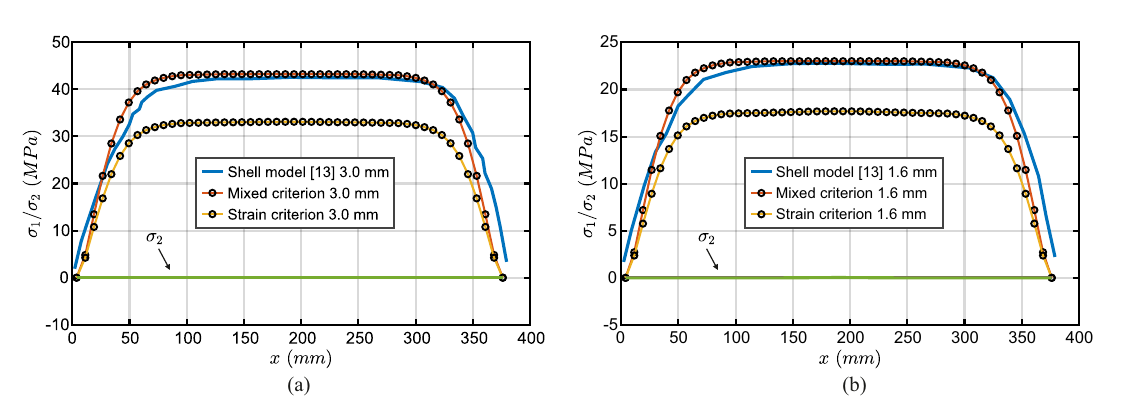}
	\caption{Principal stress comparison along the midline of a rectangular membrane under shear: (a) $u_x=\SI{3}{\milli \meter}$ and (b) $u_x=\SI{1.6}{\milli \meter}$.}
	\label{FIG:shear_sig12}
\end{figure}
As shown in Fig.~\ref{FIG:shear_sig12}, we measured and compared the principal stresses along the midline of the membrane. The results from the shell model \cite{WesleyWong2006b} serve as a reference, as they were obtained using an extremely fine shell element mesh. As observed, the mixed criterion-based model performs exceptionally well in both cases of 
$u_x = \SI{1.6}{\milli \meter}$ and $u_x = \SI{3}{\milli \meter}$, with predicted values around $\SI{23}{\mega \pascal}$ and $\SI{43}{\mega \pascal}$, respectively, which are in close agreement with the reference values. In contrast, the strain criterion-based model predicts lower principal stresses, with values around $\SI{18}{\mega \pascal}$ and $\SI{33}{\mega \pascal}$. This discrepancy may stem from the inherent limitations of the strain-based decomposition, which tends to be less precise in handling shear-dominated loading conditions. A similar phenomenon has been reported in the fracture problems~\cite{van2020strain}. As for the stress criterion-based model, due to the issues mentioned in Remark~\ref{remark1}, it fails to converge, and thus, its results are not presented. Notably, for the prediction of the second principal stress $\sigma_2$, all models yield results that are approximately zero, consistent with the shell solution. 
\begin{figure}[htbp]
	\centering
	\includegraphics[scale=0.90]{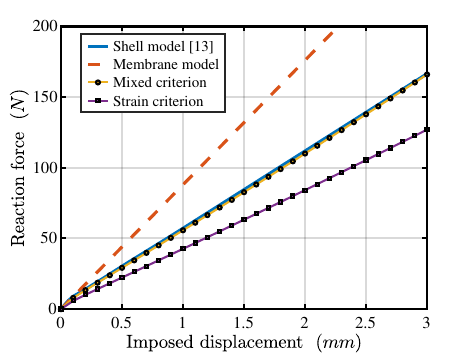}
	\caption{Shear force–displacement response.}
	\label{FIG:shear_force}
\end{figure}

Moreover, Fig.~\ref{FIG:shear_force} presents the shear force–displacement response of a rectangular membrane under shear loading, comparing different modeling approaches. The shell model serves as the reference, accurately predicting the reaction force. As observed in the figure, the slope of the force-displacement curve decreases significantly, particularly when the reaction force reaches approximately $\SI{10}{\newton}$, indicating the onset of wrinkling effects. In contrast, the membrane model exhibits a much steeper slope throughout, overestimating stiffness due to its inability to capture wrinkling effects. Among the wrinkling models, the mixed criterion closely follows the shell model, demonstrating its accuracy, while the strain criterion underestimates the reaction force. The stress criterion is not shown due to convergence issues, as previously discussed.

\subsection{A flat square membrane wrinkled by biaxial corner stretching}
\label{exp3}
In the third numerical example, we validate the proposed model by analyzing the wrinkling behavior of a square membrane subjected to corner loads $T_1$ and $T_2$. This experiment was initially conducted by \citep{WesleyWong2006}, and in their subsequent work, \citep{WesleyWong2006b} utilized a very fine shell element mesh to simulate the wrinkling phenomenon induced by biaxial stretching. Later, these shell element-based numerical results were used by many researchers as a reference to validate the accuracy of their proposed wrinkling models, as mentioned in \cite{Feng2014, Jarasjarungkiat2009}. \begin{figure}[htbp]
	\centering
	\includegraphics[scale=0.90]{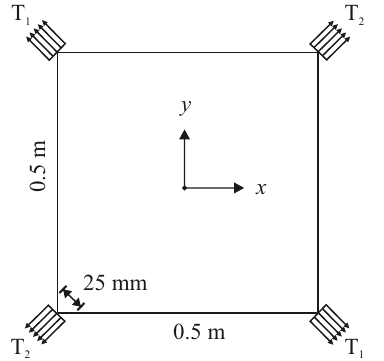}
	\caption{Flat square membrane under corner loads $T_1$ and $T_2$.}
	\label{FIG:cornerload}
\end{figure}

As shown in Fig.~\ref{FIG:cornerload}, the membrane has a side length of $\SI{0.5}{\metre}$ and a thickness of $\SI{25}{\micro\metre}$. The material properties used include Young's modulus $E = \SI{3500}{\mega\pascal}$ and Poisson's ratio $\nu = 0.31$. The membrane is subjected to equal and opposite diagonal loads $T_1$ and $T_2$ at the four corners, each applied over a width of $\SI{25}{\milli\metre}$. In the experiment, it was observed that the membrane developed wrinkles due to these stretching loads, with initial wrinkles appearing around the four corners. Subsequently, by gradually increasing $T_1$ from 5 N to 20 N while keeping $T_2$ constant, the ratio $T_1/T_2$ increased from 1 to 4, resulting in additional wrinkles forming along the diagonal direction of $T_1$.

\begin{figure}[H]
	\centering
	\includegraphics[scale=.7]{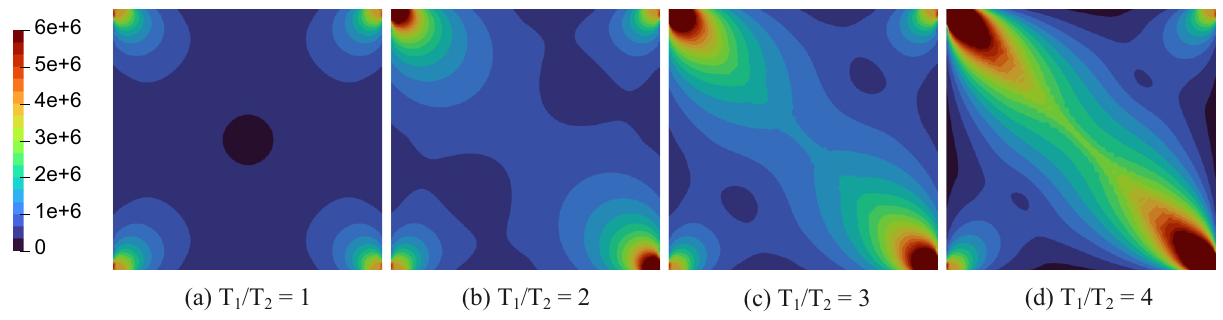}
	\caption{Contours of the first principal stresses $\sigma_1$ (\SI{}{\pascal}) predicted by mixed criterion-based model under the different values of the $T_1/T_2$-ratio.}
	\label{FIG:cornerloads_stress}
\end{figure}

In the numerical simulations, the square membrane was discretized into a $40 \times 40$ mesh of cubic membrane elements. The boundary conditions were applied as suggested by \cite{Jarasjarungkiat2009}, with the center point of the membrane fixed to prevent rigid body motion. To avoid rotation, the movement of the middle point at the top edge in the $x$-direction was constrained, while the displacement in the $z$-direction was restricted for all control points. Numerical analyses were performed for different values of the $T_1/T_2$ ratio, and the first principal stress distributions $\sigma_1$ under various load ratios $T_1/T_2$ were plotted, as shown in Figs.~\ref{FIG:cornerloads_stress}(a)-(d) for the mixed criterion-based model and in~\cite{zhang2024variationally} for the strain criterion-based model. The model based on mixed criterion appears to produce first principal stress distribution patterns comparable to those obtained from the thin shell solutions~\citep{WesleyWong2006b} and those predicted by other wrinkling models~\citep{Jarasjarungkiat2009, Feng2014}. Compared to the mixed criterion-based model, the strain criterion-based model produces similar results but shows differences at $T_1/T_2 = 2$ and $T_1/T_2 = 3$. For $T_1/T_2 = 2$, the strain criterion-based model predicts a stress distribution where the low-stress regions at the four corners become connected, whereas in the mixed criterion-based model, these regions remain separated. At $T_1/T_2 = 3$, the central stress band appears wider compared to the mixed criterion model. As for the stress criterion-based model, it converges normally for $T_1/T_2 = 1, \ 2$ and produces results similar to the mixed criterion-based model. However, for $T_1/T_2 = 3, \ 4$, convergence issues arise, and the principal stress along the diagonal approaches zero. A possible reason is that the first principal stress along the diagonal becomes negative, leading the stress criterion to classify it as a slack state.

\subsection{A square isotropic airbag wrinkled due to inflation}
\label{exp4}
In the fourth numerical example of this study, we validate the accuracy of the proposed model by simulating the three-dimensional wrinkling behavior of a square isotropic membrane under inflation pressure. This benchmark problem was first introduced by \citep{Bauer1975} and has since been widely used for validating wrinkling formulations, as documented in relevant studies~\cite{Diaby2006, Jarasjarungkiat2009, Jarasjarungkiat2008, Contri1988, kang1999finite, lee2006finite, le2021analysis, gil2007finite,zhang2020wrinkling,zhang2024variationally}.
\begin{figure}[h!]
	\centering
	\includegraphics[scale=.7]{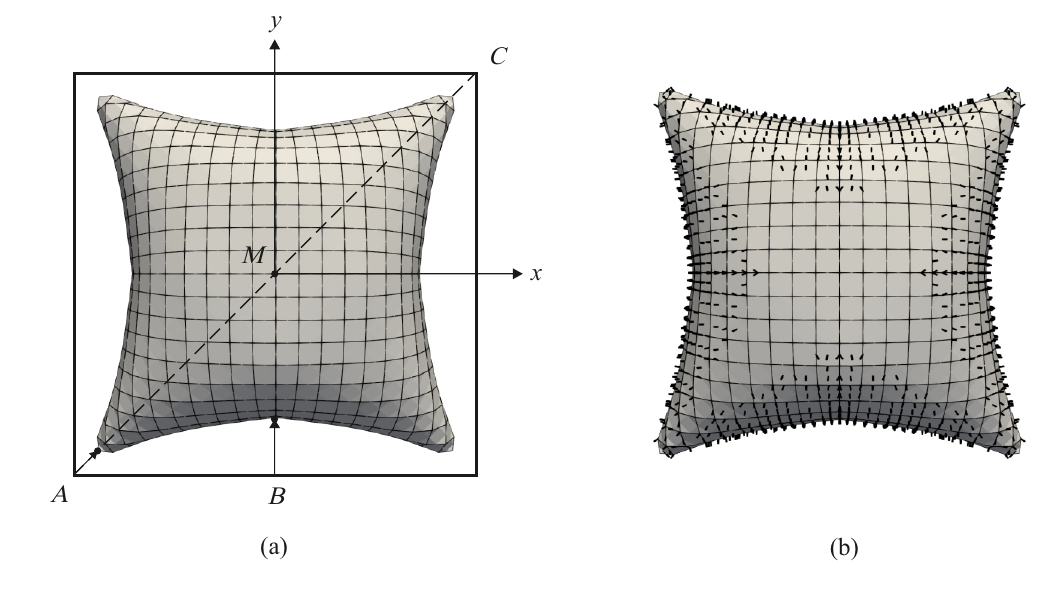}
	\caption{The  fully inflated airbag deformation (a) and its wrinkle trajectories predicted by mixed criterion-based model (b).}
	\label{FIG:airbag}
\end{figure}
As illustrated in the Fig.~\ref{FIG:airbag}(a), the initial configuration of the square airbag has a diagonal length of $AC = \SI{120}{\centi\meter}$ and a thickness of $t = \SI{0.06}{\centi\meter}$. The material properties include a Young's modulus of $E = \SI{588}{\mega\pascal}$ and a Poisson's ratio of $\nu = 0.4$. Due to symmetry, only one-quarter of the structure is simulated. Symmetric boundary conditions are applied along the inner edges, while movement along the $z$-direction at the outer edges is constrained. The airbag is subjected to a gradually increasing pressure $P$, which remains perpendicular to the airbag surface during inflation, ultimately reaching $\SI{5000}{\pascal}$. As noted in the literature, the inflation of the square airbag presents specific numerical challenges due to the singularity of the stiffness matrix at the initial stage caused by the pressure $P$ applied to the membrane surface. To address this issue, the membrane was stretched along the $x-$ and $y-$directions with death loads applied at the two outer edges. These loads were gradually reduced and removed as the pressure $P$ reached a fixed value, enabling us to avoid the initial singularity and obtain an accurate solution for the airbag inflation simulation.

For comparison with existing results, the model was discretized using bi-linear membrane elements with different mesh sizes, specifically $4 \times 4$, $5 \times 5$, $8 \times 8$, and $10 \times 10$ elements. The displacements at points $M$, $A$, and $B$, as well as the first principal stress at point $M$, were measured in each of these cases. The results obtained for each mesh size are presented in the following Tab.~\ref{tab:airbag_1} and Tab.~\ref{tab:airbag_2}:
\begin{table}[htbp]
	\centering
	\caption{Comparison of results from literature and present work for mesh refinements $4 \times 4$ and $5 \times 5$.}
	\label{tab:airbag_1}
	\resizebox{\columnwidth}{!}{%
		\begin{tabular}{lccccccccccccc}
			\toprule
			No. elements     & $4 \times 4$     &    &    &        &        &       &  & $5 \times 5$ &    &    &        &        &       \\ \cline{2-7} \cline{9-14} \\[-1em]
			& \citep{Contri1988}     & \citep{kang1999finite} & \citep{Jarasjarungkiat2009} & strain & stress & mixed &  & \citep{Contri1988} & \citep{kang1999finite} & \citep{Jarasjarungkiat2009} & strain & stress & mixed \\ 
			\midrule
			$w_M$ (m)        & 0.2090 & 0.2150   & 0.2149   & 0.2152 & 0.2159 & 0.2145  & &   0.2170  & 0.2160   & 0.2159   & 0.2163 & 0.2169 &0.2156       \\
			$r_A$ (m)        & 0.0570 &  0.0608  & 0.0972   & 0.0947 & 0.0478 &   0.0971    &  & 0.0630   &0.0594    &0.0882    & 0.0864 & 0.0479 &0.0881       \\
			$u_B$ (m)        & 0.1018 & 0.1170   & 0.1202    & 0.1199 & 0.1056 &   0.1201    & &  0.1103   & 0.1170   &0.1215    & 0.1209 & 0.1067 & 0.1213      \\
			$\sigma_M$ (MPa) & 3.4    & -   & 3.2   & 3.3      &  3.5      &    3.3   & & 3.5  &-    &3.6    & 3.6       & 3.9        & 3.6      \\ \bottomrule
		\end{tabular}%
	}
\end{table}


\begin{table}[]
	\centering
	\caption{Comparison of results from literature and present work for mesh refinements $8 \times 8$ and $10 \times 10$.}
	\label{tab:airbag_2}
	\begin{tabular}{lccccccccccc}
		\toprule
		No. elements     & $8 \times 8$     &        &        &        &        &        &  & $10 \times 10$    &        &        &        \\ \cline{2-7} \cline{9-12} 
		\\[-1em]
		& \citep{Contri1988}     & \citep{kang1999finite}     & \citep{Jarasjarungkiat2009}     & strain & stress & mixed  &  & \citep{Jarasjarungkiat2009}     & strain & stress & mixed  \\ 
		\midrule
		$w_M$ (m)        & 0.2050 & 0.2140 & 0.2166 & 0.2166 & -      & 0.2162 &  & 0.2167 & 0.2167 & -      & 0.2163 \\
		$r_A$ (m)        & 0.0470 & 0.0580 & 0.0738 & 0.0731 & -      & 0.0737 &  & 0.0692 & 0.0683 & -      & 0.0691 \\
		$u_B$ (m)        & 0.1301 & 0.1190 & 0.1227 & 0.1219 & -      & 0.1225 &  & 0.1237 & 0.1228 & -      & 0.1235 \\
		$\sigma_M$ (MPa) & 3.5    & -      & 3.8    &  3.8      & -      & 3.8    &  & 3.8    & 3.8       & -      & 3.8    \\ 
		\bottomrule
	\end{tabular}
\end{table}
The results show strong agreement with existing wrinkling models. Among them, the mixed and strain criterion-based models perform well across different mesh sizes, yielding predictions generally consistent with other models. In contrast, the stress criterion-based model tends to overestimate displacement and encounters convergence issues in certain cases. 

\begin{figure}[h!]
	\centering
	\includegraphics[scale=.9]{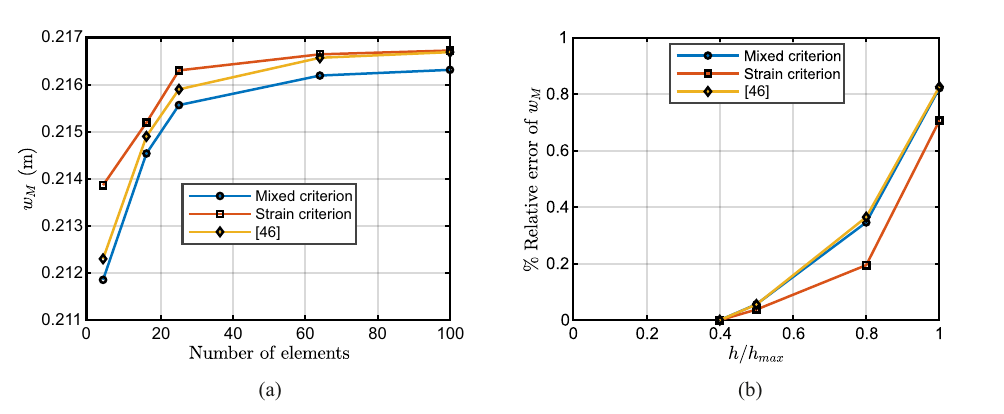}
	\caption{Convergence of vertical displacement at point M during mesh refinement (a) and its convergence rate versus relative mesh size (b).}
	\label{FIG:airbag_mesh}
\end{figure}

As shown in Fig.\ref{FIG:airbag_mesh}, we present the convergence of vertical displacement at point $M$ (Fig.\ref{FIG:airbag_mesh}(a)) and its relative error with respect to element size (Fig.\ref{FIG:airbag_mesh}(b)). The strain criterion-based model predicts values closer to the reference model \citep{Jarasjarungkiat2009} compared to the mixed criterion in terms of vertical displacement. However, when analyzing the error convergence in Fig.\ref{FIG:airbag_mesh}(b), it is evident that the mixed criterion follows a second-order ($O(h^2)$) convergence rate, consistent with the findings of reference model as reported in the literature \citep{Jarasjarungkiat2009}. The strain criterion-based model exhibits slightly larger deviations but maintains a similar convergence trend. The stress criterion-based model fails to converge and is therefore not included in the results. Furthermore, the second eigenvector of the stress tensor was used to capture the wrinkling direction of the structure, as shown in Fig.~\ref{FIG:airbag}(b), based on the mixed criterion-based model. The results indicate that wrinkles primarily form along the edges and corners, while the regions along the two diagonals remain taut without wrinkling, which is consistent with previous simulation results \citep{Diaby2006, Jarasjarungkiat2008, Jarasjarungkiat2009, Contri1988, zhang2020wrinkling}. For the strain criterion-based model, the predicted wrinkling pattern can be referenced from \citep{zhang2024variationally}.

\subsection{A square blanket wrinkled by free fall under self-weight}
\begin{figure}[h!]
	\centering
	\includegraphics[scale=0.8]{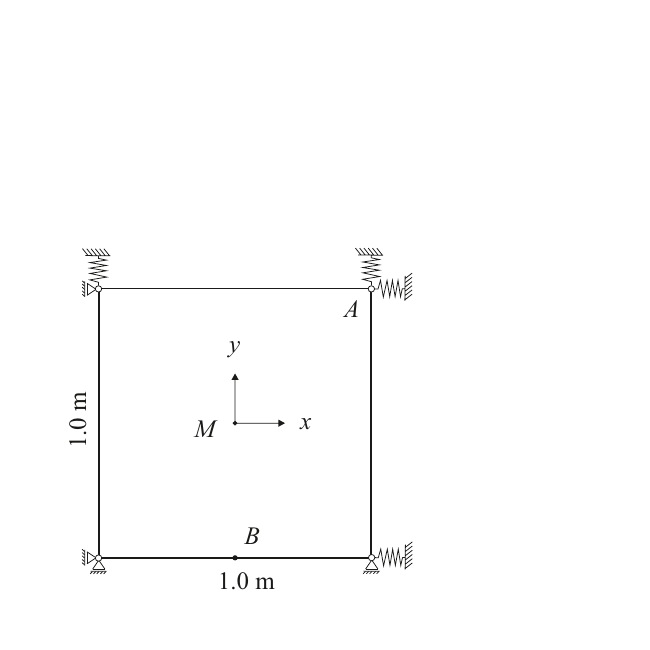}
	\caption{Hanging blanket under self-weight.}
	\label{FIG:blanket_setting}
\end{figure}

In the final example, we analyze a loose membrane subjected solely to self-weight, representing a hanging blanket supported at its corners, as depicted in Fig.\ref{FIG:blanket_setting}. The blanket has a side length of $L = \SI{1}{\meter}$ and a thickness of $t = \SI{1.77}{\milli\meter}$. The material properties include a Young’s modulus of $E = \SI{30000}{\pascal}$, Poisson’s ratio $\nu = 0.3$, and a surface density of $\rho = \SI{0.144}{\kilogram \per \milli\meter \squared}$, representing the effect of self-weight. The four corners are rigidly constrained in the $z$-direction, while elastic supports are applied in the $x$-$y$ plane, allowing significant deformations due to gravity. The elastic supports are incorporated using a penalty formulation following the approach in~\cite{Herrema2019}, with a penalty parameter of $\alpha = 10^{-2}$, equivalent to a spring stiffness of $K_{\text{spring}} = \SI{22.95}{\kilo\newton \per\meter}$. To capture the wrinkling behavior, the membrane is discretized using a $25 \times 25$ mesh of bi-quadratic elements. 

To validate the proposed models, we also perform simulations using the wrinkling model proposed in~\cite{Nakashino2005} under the same setup and mesh, and its results serve as the reference solution. Tab.~\ref{comparison_blanket} presents a comparison of the displacement components at points $A$, $B$, and $M$, along with the first principal stress at $M$, across different models. The mixed criterion-based model exhibits strong agreement with the reference solution. The strain criterion-based model predicts slightly larger displacements, particularly for $u^{M}_{z}$, whereas the stress criterion-based model also deviates but remains closer to the mixed criterion. Additionally, the strain criterion underestimates the first principal stress $\sigma_1^M$, while the stress criterion provides a closer approximation.
\begin{table}[htbp]
	\centering
	\caption{\textcolor{rv}{Comparison of the displacements and first principal stress with $\nu=0.3$}.}
	\fontsize{8pt}{8pt}\selectfont
	\setlength{\tabcolsep}{0pt} 
	\begin{tabular*}{\textwidth}{@{\extracolsep{\fill}\quad}lcccc}
		\toprule
		Mesh (25x25)           			& $u^M_z$ $(m)$  & $u^A_{x}$ $(m)$ 	& $u^B_{x}$ $(m)$  	& $\sigma^M_1$ $(MPa)$  \\
		\midrule
		Reference      					& $-0.2833$ 	& $-0.03406$ 		& $-0.01703$ 		& $642.66$ \\
		Mixed criterion               	& $-0.2833$ 	& $-0.03406$ 		& $-0.01703$ 		& $623.95$ \\
		Strain criterion               	& $-0.2956$ 	& $-0.03290$ 		& $-0.01645$ 		& $586.14$ \\
		Stress criterion               	& $-0.2887$ 	& $-0.03365$ 		& $-0.01683$ 		& $621.34$ \\
		\bottomrule
	\end{tabular*}
	\label{comparison_blanket}
\end{table}

\begin{figure}[h!]
	\centering
	\includegraphics[scale=1]{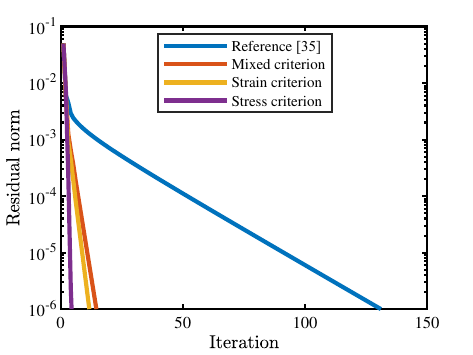}
	\caption{\textcolor{rv}{Comparison of the convergence behavior of two models during the last step with various Poisson's ratios}.}
	\label{FIG:blanket_sd_norm}
\end{figure}

Fig.~\ref{FIG:blanket_sd_norm} compares the convergence behavior of the four models. Among them, the stress criterion-based model converges the fastest, reaching numerical tolerance in $5$ iterations. The strain and mixed criterion-based models follow closely, both stabilizing within around $15$ iterations. In contrast, the reference model requires more than $100$ iterations to achieve convergence. This demonstrates that the entire series of models formulated within the spectral decomposition framework exhibit excellent convergence performance, significantly reducing computational costs compared to the reference approach.

\section{Conclusion}
\label{conclusion}
In this paper, we present a variationally consistent membrane wrinkling model based on spectral decomposition of the stress tensor, achieving notable improvements in predictive accuracy, stability, and generality. The model satisfies the uniaxial tension condition from tension field theory and aligns with the mixed wrinkling criterion. It ensures physical consistency while integrating the strengths of different decomposition approaches within a unified framework. Beyond these advancements, the model retains the superior convergence properties of a previous approach \cite{zhang2024variationally} and inherits the residual compressive stress framework, in which the negative strain energy density is scaled by a small empirical factor. Additionally, the modified strain energy density is consistently derived with respect to strain variables, leading to new stress and constitutive tensor formulations. The simplicity of the model and its adherence to variational consistency enhance its computational efficiency and robustness.

To evaluate the performance of the proposed wrinkling model, we conducted a series of analytical, numerical, and experimental benchmark tests. Overall, the proposed model delivers consistent and accurate results across various loading conditions. Additionally, we compared the proposed model with those based on other criteria. The strain criterion-based model performs well in most cases but shows reduced accuracy in shear-dominated loading scenarios. Meanwhile, the stress criterion-based model is sensitive and prone to convergence issues, making it less stable under certain conditions. These comparisons highlight the advantages of the mixed criterion-based model, particularly in terms of precision and robustness.

Future research may focus on extending the mixed criterion-based model to hyperelastic material models and anisotropic materials, further enhancing its applicability. Such extensions would expand the potential of the model to accurately capture wrinkling behavior across a broader range of materials and structural conditions.


\appendix
\section{Wrinkling model based on spectral decomposition of the strain tensor}
\label{appNew}

According to the split presented in \citep{zhang2024variationally}, the spectral decomposition is carried out on strain tensor, thereby, the strain energy density can be decomposed additionally and its positive and negative parts are defined as:
\begin{align}
\psiel\left(\mathbf{E}\right)&=\psiel^{+}\left(\mathbf{E}\right)+\psiel^{-}\left(\mathbf{E}\right), \label{psi-split} \\
\psiel^{\pm}\left(\mathbf{E}\right) &= \frac{\lambda}{2} \left(\left\langle \text{tr}\left(\mathbf{E}\right)\right\rangle ^{\pm}\right)^{2}+\mu \text{tr}\left(\left(\mathbf{E}^{\pm}\right)^{2}\right) - \frac{\lambda^2}{2\left(\lambda+2\mu\right)} \left(\left\langle \text{tr}\left(\mathbf{E}\right)\right\rangle ^{\pm}  \right)^{2}. \label{psipm}
\end{align}
Then, the stress tensor can also be split into positive and negatives parts by deriving the corresponding strain energy densities with respect to strain:
\begin{equation}
\mathbf{S}^{\pm} = (\lambda-\frac{\lambda^2}{\lambda+2\mu})\left\langle \operatorname{tr}(\mathbf{E})\right\rangle ^{\pm}\mathbf{I} + 2\mu \mathbf{E}^{\pm},  
\end{equation}
which can be rewritten as 
\begin{equation}
\label{S_strain}
\mathbf{S}^{\pm}  = \frac{E}{1-\nu^2}   \Big( \nu H^{\pm}(\operatorname{tr}(\mathbf{E}))\operatorname{tr}(\mathbf{E})\mathbf{I} + (1-\nu) H^{\pm}(E_{\alpha}) E_{\alpha}\mathbf{M}_{\alpha}\Big), 
\end{equation}
with
\begin{equation}
\lambda= \frac{E\nu }{(1+\nu )(1-2\nu )}, \quad \mu=\frac {E}{2(1+\nu )}.
\end{equation}
Since the strain tensor is symmetric, the following relation holds:
\begin{equation}
\operatorname{tr}(\mathbf{E})\mathbf{I} = (E_1+ E_2)(\mathbf{M}_1 + \mathbf{M}_2).
\end{equation}
Thus, the stress can be reformulated as
\begin{equation}
\mathbf{S}^{\pm} = \frac{E}{1-\nu^2} \bigg(  H^{\pm}(E_{\alpha}) (1-\nu) E_{\alpha}\mathbf{M}_{\alpha} + H^{\pm}(\operatorname{tr}(\mathbf{E})) \Big(\nu (E_{\alpha}\mathbf{M}_{\alpha}  + E_{\alpha}\mathbf{M}_{\beta})     \Big)     \bigg), \quad \alpha \ne \beta.
\end{equation}
Consequently, the corresponding positive and negative material tensors based on this kind of split are expressed as:
\begin{equation}
\label{C_strain}
\mathbf{\mathbb{C}}^{\pm} = \frac{E}{1-\nu^2} \bigg(     H^{\pm}(E_{\alpha}) \Big( (1-\nu)(\mathbb{Q}_{\alpha\alpha} + E_{\alpha}\mathbb{M}_{\alpha})   \Big) +  
H^{\pm}(\operatorname{tr}(\mathbf{E})) \Big(\nu (\mathbb{Q}_{\alpha \alpha } + E_{\alpha}\mathbb{M}_{\alpha}  + E_{\alpha}\mathbb{M}_{\beta} + \mathbb{Q}_{\alpha \beta})     \Big)  \bigg), \quad \alpha \ne \beta.
\end{equation}

\bibliographystyle{elsarticle-num-names} 
\bibliography{reference}

\end{document}